\documentclass[a4paper]{article}
\usepackage{amsmath}
\usepackage{amscd}
\usepackage{amstext}
\usepackage{amsfonts}
\usepackage{euscript}
\usepackage{graphicx}

\def\scr{\EuScript}

\newcommand{\C}{\mathbb{C}}
\newcommand{\ZZ}{\mathbb{Z}}
\newcommand{\NN}{\mathbb{N}}
\newcommand{\W}{\mathbb{W}}
\DeclareMathOperator{\IC}{\rm IC} \DeclareMathOperator{\Ch}{\rm Ch}
\DeclareMathOperator{\Hom}{\rm Hom}

\DeclareMathOperator{\End}{\rm End}

\DeclareMathOperator{\ann}{\rm ann}

\DeclareMathOperator{\divi}{\rm div}

\DeclareMathOperator{\im}{\rm Im}

\DeclareMathOperator{\Der}{\rm Der}

\DeclareMathOperator{\Mod}{\rm mod}

\newcommand{\derlogD}{\Der(\log D)}
\newcommand{\VCERO}{{\DX(\log D)}}
\newcommand{\VO}{{\cal V}_0}

\newcommand{\D}{{\scr D}}
\newcommand{\K}{{\scr K}}
\newcommand{\E}{{\scr E}}
\newcommand{\bE}{{\scr F}}

\newcommand{\LL}{{\scr L}}

\newcommand{\OO}{{\scr O}}

\DeclareMathOperator{\Dual}{\Bbb D}
\newcommand{\Lotimes}{\stackrel{L}{\otimes}}
\DeclareMathOperator{\DR}{DR}

\DeclareMathOperator{\adj}{adj}

\DeclareMathOperator{\Gr}{Gr}

\newcommand{\OX}{{\scr O}_X}
\newcommand{\DX}{{\scr D}_X}

\newcommand{\dx}[1]{\frac{\partial}{\partial x_{#1}}}

\newcommand{\Thetafs}{\Theta_{f,s}}

\newcommand{\CC}{\mathbb{C}}

\newcommand{\hol}{{\scr O}}

\newcounter{numero}[section]
\renewcommand{\thenumero}{(\thesection .\arabic{numero})}
\newenvironment{corolario}{\medskip
\refstepcounter{numero}\noindent {\sc  \thenumero\ Corollary.}\
\it}{\vspace{1ex}\par}

\newenvironment{teorema}{\medskip
\refstepcounter{numero}\noindent {\sc  \thenumero\ Theorem.}\
\it}{\vspace{1ex}\par}

\newenvironment{lema}{\medskip
\refstepcounter{numero}\noindent {\sc  \thenumero\ Lemma.}\
\it}{\vspace{1ex}\par}

\newenvironment{proposicion}{\medskip
\refstepcounter{numero}\noindent {\sc  \thenumero\ Proposition.}\
\it}{\vspace{1ex}\par}

\newenvironment{nota}{\medskip
\refstepcounter{numero}\noindent {\sc  \thenumero\ Remark.}\
}{\vspace{1ex}\par}

\newenvironment{ejemplo}{\medskip
\refstepcounter{numero}\noindent {\sc  \thenumero\ Example.}\
}{\vspace{1ex}\par}

\newcommand\numero{\medskip\refstepcounter{numero}\noindent{\sc \thenumero}\hspace{1em}}

\newenvironment{prueba}{
\noindent {\sc  Proof.}\ }{\hfill Q.E.D.\vspace{1ex}\par}

\newenvironment{sk-prueba}{
\noindent {\sc  Sketch of proof.}\ }{\hfill $\Box
$\vspace{1ex}\par}

\newcommand{\la}{\lambda}
\newcommand{\f}{\nu}
\newcommand{\m}{m}
\newcommand{\n}{n}

\newcommand{\ux}{\underline{x}}
\newcommand{\umu}{\underline{\mu}}

\newcommand{\ub}{\underline{b}}
\newcommand{\om}{\omega}

\newcommand{\comb}[2]{\left(\!\!\!\!
                      \begin{array}{c}
                       #1 \\
                       #2  \\
                      \end{array} \!\!\!\!\right)}

\title{Algebraic computation of some intersection D-modules}

\author{F. J. Calder\'{o}n Moreno and L. Narv\'{a}ez Macarro\thanks{
The authors are partially supported by MTM2004-07203-C02-01
 and FEDER.}}

\date{March, 2006}

\begin{document}

\maketitle

\begin{abstract}
Let $X$ be a complex analytic manifold, $D\subset X$ a locally
quasi-homo\-ge\-neous free divisor, $\E$ an integrable logarithmic
connection with respect to $D$ and $\LL$ the local system of the
horizontal sections of $\E$ on $X-D$. In this paper we give an
algebraic description in terms of $\E$ of the regular holonomic
$\DX$-module whose de Rham complex is the intersection complex
associated with $\LL$. As an application, we perform some effective
computations in the case of quasi-homogeneous plane curves.
\end{abstract}

\section*{Introduction}

On a complex analytic manifold, intersection complexes associated
with irreducible local systems on a dense open regular subset of a
closed analytic subspace are the simple pieces which form any
perverse sheaf. The Riemann-Hilbert correspondence allows us to
consider the regular holonomic D-modules which correspond to these
intersection complexes, that we call ``intersection D-modules". They
are the simple pieces which form any regular holonomic D-module.
Whereas intersection complexes are topological objects, intersection
D-modules are algebraic: they are given by a system of partial linear
differential equations with holomorphic coefficients.

Intersection complexes can be constructed by an important operation:
the intermediate direct image. Its description in terms of Verdier
duality and usual derived direct images can be algebraically
interpreted in the category of holonomic regular D-modules by using
the deep properties of the de Rham functor. We need to compute
localizations and D-duals.

This can be effectively done, in principle, by using the general
available algorithms in \cite{oaku_dk_97,MR1769663,MR1808827}, but
in the case of integrable logarithmic connections along a locally
quasi-homogeneous free divisor, we exploit the logarithmic point of
view
\cite{calde_ens,calde_nar_compo,calde_nar_fourier,cas_ucha_stek,cas_ucha_exper,torre-45-bis,torre-overview}
to previously obtain a \underline{general algebraic description} of
their associated intersection D-modules, from which we can easily
derive effective computations.

The main ingredients we use are the duality theorem proved in
\cite{calde_nar_fourier} and the logarithmic comparison theorem for
arbitrary integrable logarithmic connections proved in
\cite{calde_nar_LCTILC}, both with respect to locally
quasi-homogeneous free divisors.

The algorithmic treatment of the computations in this paper will be
developed elsewhere.

Let us now comment on the content of this paper.

In section \ref{sec:1} we remind the reader of the basic notions and
notations and we review our previous results on logarithmic
$\D$-modules with respect to free divisors. We recall the
logarithmic comparison theorem for arbitrary integrable logarithmic
connections from \cite{calde_nar_LCTILC}, and we give the theorem
describing the intersection D-module associated with an integrable
logarithmic connection along a locally quasi-homogeneous free
divisor.

In section \ref{sec:2}, given a locally quasi-homogeneous free
divisor $D$ with a reduced local equation $f=0$ and a cyclic
integrable logarithmic connection $\E$ with respect to $D$, we
explicitly describe a presentation of $\D[s]\cdot \left(\E
f^s\right)$ over $\D[s]$ in terms of a presentation of $\E$ over the
ring of logarithmic differential operators. This description will be
useful in order to compute the Bernstein-Sato polynomials associated
with $\E$.

In section \ref{sec:3}, the general results of the previous section
are explicitly written down in the case of a family of integrable
logarithmic connections with respect to a quasi-homogeneous plane
curves.

In section \ref{sec:4} we perform some explicit computations with
respect to a cusp.

We wish to thank H\'el\`ene Esnault who, because of a question about
our paper \cite{calde_nar_fourier}, drew our attention to computing
intersection D-modules. We also thank Tristan Torrelli for helpful
information about the Bernstein-Sato functional equations and for
some comments on a previous version of this paper.

\section{Logarithmic connections with respect to a free divisor: theoretical set-up}\label{loco}
\label{sec:1} Let $X$ be a $n$-dimensional complex analytic manifold
and $D\subset X$ a hypersurface, and let us denote by $j: U=X-D
\hookrightarrow X$ the corresponding open inclusion.

 We say that $D$ is a {\em free divisor} \cite{ksaito_log} if the
$\OX$-module $\Der(\log D)$ of logarithmic vector fields with respect
to $D$ is locally free (of rank $n$), or equivalently if the
$\OX$-module $\Omega^1_X(\log D)$ of logarithmic 1-forms with respect
to $D$ is locally free (of rank $n$).

Normal crossing divisors, plane curves, free hyperplane arrangements
(e.g. the union of reflecting hyperplanes of a complex reflection
group), discriminant of stable mappings or bifurcation sets are
examples of free divisors.

We say that $D$ is quasi-homogeneous at $p\in D$ if there is a
system of local coordinates  $\ux$ centered at $p$ such that the
germ $(D,p)$ has a reduced weighted homogeneous defining equation
(with strictly positive weights) with respect to $\ux$. We say that
$D$ is locally quasi-homogeneous if it is so at each point $p\in D$.

Let us denote by $\DX(\log D)$ the $0$-term of the
Malgrange-Kashiwara filtration with respect to $D$ on the sheaf $\DX$
of linear differential operators on $X$. When $D$ is a free divisor,
the first author has proved in \cite{calde_ens} that $\DX(\log D)$ is
the universal enveloping algebra of the Lie algebroid $\Der(\log D)$,
and then it is coherent and has noetherian stalks of finite global
homological dimension. Locally, if $\{\delta_1,\dots,\delta_n\}$ is a
local basis of the logarithmic vector fields on an open set $V$, any
differential operator in $\Gamma(V,\DX(\log D))$ can be written in a
unique way as a finite sum
$$ \sum_{\substack{\alpha \in\NN^n\\ |\alpha|\leq d}} a_{\alpha}
\delta_1^{\alpha_1}\cdots \delta_n^{\alpha_n},$$where the
$a_{\alpha}$ are holomorphic functions on $V$.

From now on, let us assume that $D$ is a free divisor.

We say that $D$ is a {\em Koszul free} divisor \cite{calde_ens} at a
point $p\in D$ if the symbols of any (some) local basis
$\{\delta_1,\dots,\delta_n\}$ of $\Der(\log D)_p$ form a regular
sequence in $\Gr {\scr D}_{X,p}$. We say that $D$ is a {\em Koszul
free} divisor if it is so at any point $p\in D$. Actually, as M.
Schulze pointed out, Koszul freeness is equivalent to holonomicity in
the sense of \cite{ksaito_log}.

Plane curves and locally quasi-homogeneous free divisors (e.g. free
hyperplane arrangements or discriminant of stable mappings in
Mather's ``nice dimensions") are example of Koszul free divisors
\cite{calde_nar_LQHKF}.

A {\em logarithmic connection} with respect to $D$ is a locally free
$\OX$-module $\E$ endowed with:\\ -) a $\C$-linear morphism
(connection) $ \nabla': \E\xrightarrow{}
\E\otimes_{\OX}\Omega^1_X(\log D)$, satisfying $\nabla'(ae) = a
\nabla'(e) + e\otimes da$, for any section $a$
of $\OX$ and any section $e$ of $\E$,\\
or equivalently, with\\
-) a left $\OX$-linear morphism $\nabla:\derlogD \xrightarrow{}
\End_{\C_X}(\E)$ satisfying the Leibniz rule $\nabla(\delta)(ae) =
a\nabla(\delta)(e) + \delta(a)e$, for any logarithmic vector field
$\delta$, any section $a$ of $\OX$ and any section $e$ of $\E$.

The integrability of $\nabla'$ is equivalent to the fact that
$\nabla$ preserve Lie brackets. Then, we know from \cite{calde_ens}
that giving an integrable logarithmic connection on a locally free
$\OX$-module $\E$ is equivalent to extending its original
$\OX$-module structure to a left $\DX(\log D)$-module structure, and
so integrable logarithmic connections are the same as left $\DX(\log
D)$-modules which are locally free of finite rank over $\OX$.

Let us denote by $\OX(\star D)$ the sheaf of meromorphic functions
with poles along $D$. It is a holonomic left $\DX$-module.

The first examples of integrable logarithmic connections (ILC for
short) are the invertible $\OX$-modules $\OX(mD)\subset \OX(\star
D)$, $m\in\ZZ$, formed by the meromorphic functions $h$ such that
$\divi (h) + mD \geq 0$.

If $f=0$ is a reduced local equation of $D$ at $p\in D$ and
$\delta_1,\dots,\delta_n$ is a local basis of $\Der(\log D)_p$ with
$\delta_i(f)=\alpha_i f$, then $f^{-m}$ is a local basis of
$\OO_{X,p}(mD)$ over $\OO_{X,p}$ and we have the following local
presentation over $\D_{X,p}(\log D)$ (\cite{calde_ens}, th. 2.1.4)
\begin{equation}\label{eq:presenta}
\OO_{X,p}(mD) \simeq \D_{X,p}(\log D)/\D_{X,p}(\log
D)(\delta_1+m\alpha_1,\dots, \delta_n+m\alpha_n).\end{equation}

\numero \label{nume:1} For any ILC $\E$ and any integer $m$, the
locally free $\OX$-modules $\E(mD):= \E \otimes_{\OX} \OX(mD)$ and
$\E^* := \Hom_{\OX}(\E,\OX)$ are endowed with a natural structure of
left $\DX(\log D)$-module, where the action of logarithmic vector
fields is given by
\begin{equation} \label{eq:oper-ext}
(\delta h)(e) = - h(\delta e) + \delta (h(e)), \quad \delta \left(
e\otimes a\right) = (\delta e)\otimes a + e\otimes \delta (a)
\end{equation}
for any logarithmic vector field $\delta$, any local section $h$ of
$\Hom_{\OX}(\E,\OX)$, any local section $e$ of $\E$ and any local
section $a$ of $\OX(mD)$ (cf. \cite{calde_nar_fourier}, \S 2). Then
$\E(mD)$ and $\E^*$ are ILC again, and the usual isomorphisms
$$ \E(mD)(m'D) \simeq \E((m+m')D),\quad \E(mD)^* \simeq \E^*(-mD)$$are
$\DX(\log D)$-linear.

\numero \label{num:0} If $D$ is Koszul free and $\E$ is an ILC, then
the complex $\DX\Lotimes_{\DX(\log D)} \E$ is concentrated in degree
$0$ and its $0$-cohomology $\DX\otimes_{\DX(\log D)} \E$ is a
holonomic $\DX$-module (see \cite{calde_nar_fourier}, prop. 1.2.3).

If $\E$ is an ILC, then $\E(\star D)$ is a meromorphic connection
(locally free of finite rank over $\OX(\star D)$) and then it is a
holonomic $\DX$-module (cf. \cite{meb_nar_dmw}, th. 4.1.3).
Actually, $\E(\star D)$ has regular singularities on the smooth part
of $D$ (it has logarithmic poles! \cite{del_70}) and then it is
regular everywhere \cite{meb-cimpa-2}, cor. 4.3-14, which means that
if $\LL$ is the local system of horizontal sections of $\E$ on
$U=X-D$, the canonical morphism
$$ \Omega_X^{\bullet}(\E(\star D)) \xrightarrow{} R j_* \LL$$is an
isomorphism in the derived category.

For any ILC $\E$, or even for any left $\DX(\log D)$-module (without
any finiteness property over $\OX$), one can define its logarithmic
de Rham complex $\Omega^{\bullet}_X(\log D)(\E)$ in the classical
way (cf. \cite[def.~I.2.15]{del_70}), which is a subcomplex of
$\Omega^{\bullet}_X(\E(\star D))$. It is clear that both complexes
coincide on $U$.

For any ILC $\E$ and any integer $m$, $\E(mD)$ is a sub-$\DX(\log
D)$-module of the regular holonomic $\DX$-module $\E(\star D)$, and
then we have a canonical morphism in the derived category of left
$\DX$-modules
$$ \rho_{\E,m}: \DX \Lotimes_{\DX(\log D)} \E(mD)\to \E(\star D),$$
given by $\rho_{\E,m}(P\otimes e') = Pe'$.

Since $\E(m'D)(mD) = \E((m+m')D)$ and $\E(m'D)(\star D) = \E(\star
D)$, we can identify morphisms $\rho_{\E(m'D),m}$ and
$\rho_{\E,m+m'}$.

For any bounded complex $\K$ of sheaves of $\C$-vector spaces on $X$,
let us denote by $\K^{\vee}= R \Hom_{\C_X}(\K,\C_X)$ its Verdier
dual.

The dual local system $\LL^{\vee}$ appears as the local system of
the horizontal sections of the dual ILC $\E^*$.

We have the following theorem (see \cite[th.~4.1]{calde_nar_fourier}
and \cite[th.~(2.1.1)]{calde_nar_LCTILC}):

\begin{teorema} \label{teo:crit-LCT} Let $\E$ be an ILC (with respect to the divisor $D$)
and let $\LL$ be the local system of its horizontal sections on
$U=X-D$. The following properties are equivalent:
\begin{enumerate}
\item[1)] The canonical morphism $\Omega_X^{\bullet}(\log D)(\E)
\to R j_* \LL$ is an isomorphism in the derived category of
complexes of sheaves of complex vector spaces.
\item[2)] The inclusion $\Omega_X^{\bullet}(\log D)(\E) \hookrightarrow
\Omega_X^{\bullet}(\E(\star D))$ is a quasi-isomorphism.
\item[3)] The morphism $\rho_{\E,1}: \DX\Lotimes_{\DX(\log D)} \E(D)\to \E(\star D)$
is an isomorphism in the derived category of left $\DX$-modules.
\item[4)] The complex $\DX\Lotimes_{\DX(\log D)} \E(D)$ is
concentrated in degree $0$ and the $\DX$-module
$\DX\otimes_{\DX(\log D)} \E(D)$ is holonomic and isomorphic to its
localization along $D$.
\end{enumerate}
Moreover, if $D$ is a Koszul free divisor, the preceding properties
are also equivalent to:
\begin{enumerate}
\item[5)] The canonical morphism $j_!\LL^{\vee} \to
\Omega_X^{\bullet}(\log D)(\E^*(-D))$ is an isomorphism in the
derived category of complexes of sheaves of complex vector spaces.
\end{enumerate}
\end{teorema}

For $D$ a locally quasi-homogeneous free divisor and $\E=\OX$, the
equivalent properties in theorem \ref{teo:crit-LCT} hold: this is the
so called ``logarithmic comparison theorem" \cite{cas_mond_nar_96}
(see also \cite[th.~4.4]{calde_nar_fourier} and
\cite[cor.~(2.1.3)]{calde_nar_LCTILC} for other proofs based on
D-module theory).
\medskip

\numero \label{nume:BS} Let $\E$ be an ILC (with respect to $D$) and
$p$ a point in $D$. Let $f\in\hol=\hol_{X,p}$ be a reduced local
equation of $D$ and let us write $\D= \D_{X,p}$, $\VO=\DX(\log D)_p$
and $E=\E_p$. We know from \cite[lemma~(3.2.1)]{calde_nar_LCTILC}
that the ideal of polynomials $b(s)\in\C[s]$ such that
$$ b(s) Ef^s \subset \D[s]\cdot \left(Ef^{s+1}\right) \left( \subset E[f^{-1},s]f^s \right)$$
is generated by a non constant polynomial $b_{\E,p}(s)$. By the
coherence of the involved objects we deduce that $b_{\E,q}(s)\ |\
b_{\E,p}(s)$ for $q\in D$ close to $p$.

If $b_{\E,p}(s)$ has some integer root, let us call $\kappa(\E,p)$
the minimum of those roots. If not, let us write $\kappa(\E,p)=
+\infty$.

Let us call
$$ \kappa(\E) = \inf \{\kappa(\E,p)\ |\ p\in D\} \in \ZZ \cup \{\pm \infty\}.$$

From now on let us suppose that $D$ is a locally quasi-homogeneous
free divisor.

\begin{teorema}  \label{teo:main} Under the above hypothesis, if
$\kappa(\E) > -\infty$, then the morphism
\begin{equation} \label{eq:1}
\rho_{\E,k}: \DX \Lotimes_{\DX(\log D)} \E(kD)\xrightarrow{}
\E(\star D)
\end{equation}
is an isomorphism in the derived category of left $\DX$-modules, for
all $k\geq -\kappa(\E)$.
\end{teorema}

\begin{prueba} It is a straightforward consequence of
\cite{calde_nar_LQHKF}, \cite[th.~5.6]{calde_nar_compo} and theorem
(3.2.6) of \cite{calde_nar_LCTILC} and its proof.
\end{prueba}

Let us note that the hypothesis $\kappa(\E) > -\infty$ in theorem
\ref{teo:main} holds locally on $X$.
\medskip

In the situation of theorem \ref{teo:main}, if $\LL$ is the local
system of the horizontal sections of $\E$ on $U=X-D$, then the
derived direct image $R j_* \LL$ is canonically isomorphic (in the
derived category) to the de Rham complex of the holonomic
$\DX$-module $\DX \otimes_{\DX(\log D)} \E(kD)$:
\begin{eqnarray*}
& \DR \left(\DX \otimes_{\DX(\log D)} \E(kD)\right) = \DR
\left(\DX \Lotimes_{\DX(\log D)} \E(kD)\right) \simeq&\\
& \DR \E(\star D) \simeq \Omega_X^{\bullet}(\E(\star D)) \simeq R
j_* \LL.&
\end{eqnarray*}

Proceeding as above for the dual ILC $\E^*$, we find that if
$\kappa(\E^*)> -\infty$, then we have that the canonical morphism
$$ \DR \left(\DX \otimes_{\DX(\log D)} \E^*(k'D)\right) \xrightarrow{}
R j_* \LL^{\vee}$$ is an isomorphism in the derived category for
$k'\geq -\kappa(\E^*)$.

Let us denote by
\begin{equation} \label{eq:2} \varrho_{\E,k,k'}:\DX \otimes_{\DX(\log D)} \E((1-k')D)
\xrightarrow{} \DX \otimes_{\DX(\log D)} \E(kD),
\end{equation}
the $\DX$-linear morphism induced by the inclusion $\E((1-k')D)
\subset \E(kD)$, $1-k'\leq k$, and by $\IC_X(\LL)$ the intersection
complex of Deligne-Goresky-MacPherson associated with $\LL$, which
is described as the intermediate direct image $j_{!*}\LL$, i.e. the
image of $j_!\LL \to R j_*\LL$ in the category of perverse sheaves
(cf. \cite{bbd_83}, def. 1.4.22).

The following theorem describes the ``intersection $\DX$-module"
corresponding to $\IC_X(\LL)$ by the Riemann-Hilbert correspondence
of Mebkhout-Kashiwara \cite{kas_RH,meb_I_84,meb_II_84}.

\begin{teorema} \label{teo:2} Under the above hypothesis,
we have a canonical isomorphism in the category of perverse sheaves
on $X$,
$$ \IC_X(\LL) \simeq \DR \left( \im \varrho_{\E,k,k'}\right),$$
for $k\geq -\kappa(\E)$, $k'\geq -\kappa(\E^*)$ and $1-k'\leq k$.
\end{teorema}

\begin{prueba} Using our duality results in \cite[\S 3]{calde_nar_fourier},
the Local Duality Theorem for holonomic $\DX$-modules
(\cite{meb_formalisme}, ch. I, th. (4.3.1); see also \cite{nar-ldt})
and theorem \ref{teo:main}, we obtain
\begin{eqnarray*}
& \DR \left(\DX \otimes_{\DX(\log D)} \E((1-k')D)\right) \simeq \DR
 \left(\DX \otimes_{\DX(\log D)} \E^*(k'D)^*(D) \right) \simeq &\\
& \DR \left( \Dual_{\DX} \left(\DX \otimes_{\DX(\log D)} \E^*(k'D)
\right) \right) \simeq \left[\DR \left( \DX \otimes_{\DX(\log D)}
\E^*(k'D) \right)\right]^{\vee} \simeq&\\ & \left[R j_*
\LL^{\vee}\right]^{\vee} \simeq j_! \LL.&
\end{eqnarray*}
On the other hand, the canonical morphism $j_! \LL \to R j_* \LL$
corresponds, through the de Rham functor, to the $\DX$-linear
morphism $\varrho_{\E,k,k'}$, and the theorem is a consequence of
the Riemann-Hilbert correspondence which says that the de Rham
functor establishes an equivalence of abelian categories between the
category of regular holonomic $\DX$-modules and the category of
perverse sheaves on $X$.
\end{prueba}

\begin{nota} For $\E=\OX$, one has $\E^*=\OX$ and there are examples where morphisms
$\rho_{\OX,k}$  in (\ref{eq:1}) are never isomorphisms
(\cite{calde_nar_fourier}, ex. 5.3). Nevertheless, for $k=k'=1$ the
image of the morphism
$$ \varrho_{\OX,1,1}:\DX \otimes_{\DX(\log D)} \OX
\xrightarrow{} \DX \otimes_{\DX(\log D)} \OX(D)$$ is always
(canonically isomorphic to) $\OX$, which is the regular holonomic
$\DX$-module corresponding by the Riemann-Hilbert correspondence to
$\IC_X(\CC_U) = \CC_X$, where $\CC_U$ is the local system of
horizontal sections of $\OX$ on $U$. To see this, let us work locally
as in (\ref{eq:presenta}). Then, morphism $\varrho_{\OX,1,1}$ is
given at point $p$ by
$$ \overline{P}\in \D_{X,p}/\D_{X,p}(\delta_1,\dots, \delta_n) \mapsto \overline{Pf}\in \D_{X,p}
/\D_{X,p}(\delta_1+\alpha_1,\dots, \delta_n+\alpha_n)$$and the stalk
at $p$ of $\im \varrho_{\OX,1,1}$ is given by $\D_{X,p}/J$ where $J$
is the left ideal
$$ J = \{P\in \D_{X,p}\ |\ Pf \in \D_{X,p}(\delta_1+\alpha_1,\dots,
\delta_n+\alpha_n)\}.$$By Saito's criterion \cite{ksaito_log} we can
suppose
$$ \left(\begin{array}{c} \delta_1\\
\vdots \\ \delta_n \end{array}\right) = A
\left(\begin{array}{c} \frac{\partial}{\partial x_1}\\ \vdots \\
\frac{\partial}{\partial x_n} \end{array}\right)
$$
where $A$ is a $n\times n$ matrix with entries in $\OO_{X,p}$ and
$\det A = f$. Writing  $B= \adj (A)^t$ we obtain
$$ B \left(\begin{array}{c} \delta_1\\
\vdots \\ \delta_n \end{array}\right) = f \left(\begin{array}{c} \frac{\partial}{\partial x_1}\\ \vdots \\
\frac{\partial}{\partial x_n} \end{array}\right)\quad
\stackrel{\text{eval. on $f$}}{\rightsquigarrow}\quad
B \left(\begin{array}{c} \alpha_1\\
\vdots \\ \alpha_n \end{array}\right) = \left(\begin{array}{c} \frac{\partial f}{\partial x_1}\\ \vdots \\
\frac{\partial f}{\partial x_n} \end{array}\right).
$$
 Then
$$ \left(\begin{array}{c} \frac{\partial}{\partial x_1}\\ \vdots \\
\frac{\partial}{\partial x_n} \end{array}\right) f = f \left(\begin{array}{c} \frac{\partial}{\partial x_1}\\ \vdots \\
\frac{\partial}{\partial x_n} \end{array}\right) + \left(\begin{array}{c} \frac{\partial f}{\partial x_1}\\ \vdots \\
\frac{\partial f}{\partial x_n} \end{array}\right) = \cdots = B \left(\begin{array}{c} \delta_1+\alpha_1\\
\vdots \\ \delta_n + \alpha_n\end{array}\right)$$ and
$\frac{\partial}{\partial x_i}\in J$ for $i=1,\dots,n$. Since $J$ is
is not the total ideal, we deduce by maximality that $J$ is the ideal
generated by the $\frac{\partial}{\partial x_i}$ and
$\D_{X,p}/J\simeq \OO_{X,p}$. To conclude, one easily sees, from the
fact that morphism $\varrho_{\OX,1,1}$ factors through
$$a\in \OX \mapsto 1\otimes a \in \DX\otimes_{\VCERO} \OX(D)$$
[it is $\DX$-linear since, for any derivation $\delta$ and any
holomorphic function $a$, $\delta(1\otimes a) = \delta\otimes a =
\delta \otimes (ff^{-1}a) = (\delta f)\otimes (f^{-1}a) = 1\otimes
(\delta f)(f^{-1} a) = 1\otimes (\delta a)$] that the isomorphisms
above at different $p$ glue together and give a global isomorphism
$\im \varrho_{\OX,1,1} \simeq \OX$.
\medskip

This example suggests studying the comparison between $\DR(\im
\varrho_{\E,k,k'})$, $k,k'\gg 0$, and $\IC_X(\LL)$ in theorem
\ref{teo:2}, independent of the fact that $\rho_{\E,k}$ and
$\rho_{\E^*,k'}$ are isomorphisms or not.
\end{nota}

\section{Bernstein-Sato polynomials for cyclic integrable logarithmic
connections}  \label{sec:2}

In the situation of \ref{nume:BS}, let us assume that $E$ is a cyclic
$\VO$-module generated by an element $e\in E$. The following result
is proved in \cite[prop.~(3.2.3)]{calde_nar_LCTILC}.

\begin{proposicion} \label{prop:BS} Under the above conditions, the polynomial $b_{\E,p}(s)$
coincides with the Bernstein-Sato polynomial $b_e(s)$ of $e$ with
respect to $f$, where $e$ is considered to be an element of the
holonomic $\D$-module $E[f^{-1}]$ (cf. \cite{kas_78}).
\end{proposicion}

\numero \label{nume:thetafs} Let $\Thetafs \subset \D[s]$ be the set
of operators in $\ann_{\D[s]} f^s$ of total order (in $s$ and in the
derivatives) $\leq 1$. The elements of $\Thetafs$ are of the form
$\delta - \alpha s$ with $\delta\in \Der_{\C}(\hol)$, $\alpha\in
\hol$ and $\delta(f)=\alpha f$. In particular $\Thetafs \subset
\VO[s]$.
\smallskip

The $\hol$-linear map
$$\textstyle \delta \in \Der(\log D)_p \mapsto \delta - \frac{\delta(f)}{f}
s\in \Thetafs$$ is an isomorphism of Lie-Rinehart algebras over
$(\C,\hol)$ and extends to a unique ring isomorphism $\Phi: \VO[s]
\xrightarrow{} \VO[s]$ with $\Phi(s)=s$ and $\Phi(a)=a$ for all
$a\in\hol$. Let us note that $\Phi^{-1}(\delta) = \delta +
\frac{\delta(f)}{f} s$ for each $\delta\in \Der(\log D)_p$.
\smallskip

It is clear that $E[s]f^s$ is a sub-$\VO[s]$-module of
$E[s,f^{-1}]f^s$ and that for any $P\in \VO[s]$ and any $e'\in E[s]$,
the following relation holds
\begin{equation} \label{eq:rel-main}
 (P e') f^s = \Phi(P)(e' f^s).
\end{equation}

\begin{proposicion} \label{prop:annVOefs} Under the above conditions, the following
relation holds
$$ \ann_{\VO[s]}\left( e f^s\right) = \VO[s]\cdot \Phi\left( \ann_{\VO}
e\right).$$
\end{proposicion}

\begin{prueba} The inclusion $\supset$ comes from
(\ref{eq:rel-main}). For the other inclusion, let $Q\in
\ann_{\VO[s]}\left( e f^s\right)$ and let us write $\Phi^{-1}(Q)=
\sum_{i=1}^d P_i s^i$ with $P_i\in \VO$. We have
$$ 0 = Q\left( e f^s \right) = \left( \Phi^{-1}(Q) e\right) f^s =
\left( \sum_{i=1}^d (P_i e) s^i \right) f^s$$ and then $P_i\in
\ann_{\VO} e$. Therefore
$$ Q =\Phi \left( \sum_{i=1}^d P_i s^i \right) = \sum_{i=1}^d
\Phi(P_i) s^i \in \VO[s]\cdot \Phi\left( \ann_{\VO} e\right).$$
\end{prueba}

\begin{proposicion} \label{prop:annDefs} Under the above conditions, if $D$
is a locally quasi-homogeneous
free divisor, then $$\ann_{\D[s]} (e f^s) = \D[s]\cdot \ann_{\VO[s]}
(e f^s).$$
\end{proposicion}

\begin{prueba} From (\ref{eq:rel-main}) we know that $E[s]f^s =
\VO[s]\cdot (e f^s)$, and from \cite[cor.~(3.1.2)]{calde_nar_LCTILC}
we know that the morphism
$$\rho_{E,s}: P\otimes (e'f^s) \in \D[s]\otimes_{\VO[s]} E[s]f^s
\mapsto P(e' f^s) \in \D[s]\cdot (E[s] f^s) = \D[s]\cdot (e f^s)$$
is an isomorphism of left $\D[s]$-modules. Therefore
$$\ann_{\D[s]} (e f^s) = \D[s]\cdot \ann_{\VO[s]}
(e f^s).$$
\end{prueba}

\begin{corolario} \label{cor:annDefs} Under the above conditions, if $D$ is a
locally quasi-homo\-ge\-neous free divisor, then
$$\ann_{\D[s]} (e f^s) = \D[s]\cdot \Phi\left( \ann_{\VO}
e\right).$$
\end{corolario}

\begin{prueba} It follows from propositions \ref{prop:annVOefs} and
\ref{prop:annDefs}.
\end{prueba}

\begin{nota} \label{nota:a}
Theorems \ref{teo:main} and \ref{teo:2}, proposition \ref{prop:annDefs}
and corollary \ref{cor:annDefs} remain true if we only assume that
our divisor $D$ is Koszul free and of commutative linear type, i.e.
its jacobian ideal is of linear type (see \cite[\S
3]{calde_nar_LCTILC}).
\end{nota}

\begin{nota}
As we shall see in sections \ref{sec:3} and \ref{sec:4}, theorem
\ref{teo:2}, proposition \ref{prop:BS} and corollary
\ref{cor:annDefs} provide an effective method of computing the
intersection $\DX$-module corresponding to $\IC_X(\LL)$ in terms of
the ILC $\E$, at least if $D$ is a locally quasi-homogeneous free
divisor, or more generally, if $D$ is Koszul free and of commutative
linear type (see remark \ref{nota:a}).
\end{nota}

\begin{nota} In the particular case of $\E =\OX$ and $E=\hol$,
corollary \ref{cor:annDefs} says that
$$ \ann_{\D[s]} (f^s) = \D[s]\cdot \left(\delta_1 - \alpha_1
s,\dots,\delta_n -\alpha_n s\right),$$ where
$\delta_1,\dots,\delta_n$ is a local basis of $\Der(\log D)_p$ and
$\delta_i(f)=\alpha_i f$ (see corollary 5.8, (b) in
\cite{calde_nar_compo}).
\end{nota}

\begin{ejemplo} Let us suppose that $D\subset X$ is a
non-necessarily free divisor and let $f=0$ be a reduced local
equation of $D$ at a point $p\in D$. Let
$\{\delta_1,\dots,\delta_m\}$ a system of generators of $\Der(\log
D)_p$ and let us write $\delta_i(f)=\alpha_i f$.

Let us call $\ann_{\D[s]}^{(1)}(f^s)$ the ideal of $\D[s]$ generated
by $\Thetafs$ (see \ref{nume:thetafs}):
$$ \ann_{\D[s]}^{(1)}(f^s) = \D[s]\cdot\left(\delta_1 - \alpha_1
s,\dots,\delta_m -\alpha_m s\right)\subset \ann_{\D[s]}(f^s).$$

The Bernstein functional equation for $f$
$$ b(s) f^s = P(s) f^{s+1}$$ means that the operator $b(s) - P(s) f$
belongs to the annihilator of $f^s$ over $\D[s]$. Then, an explicit
knowledge of the ideal $\ann_{\D[s]}(f^s)$ allows us to find $b(s)$
by computing the ideal $$\C[s]\cap \left( \D[s]\cdot f +
\ann_{\D[s]}(f^s) \right),$$ (see \cite{oaku_dk_97}). However, the
ideal $\ann_{\D[s]}(f^s)$ is in general difficult to compute.

When $D$ is a locally quasi-homogeneous free divisor, or more
generally, a divisor of differential linear type
(\cite{calde_nar_LCTILC}, def. (1.4.5) ), $\ann_{\D[s]}(f^s)=
\ann_{\D[s]}^{(1)}(f^s)$ and the computation of $b(s)$ is in
principle easier.

But there are other examples where the Bernstein polynomial $b(s)$
belongs to
$$\C[s]\cap \left( \D[s]\cdot f +
\ann_{\D[s]}^{(1)}(f^s) \right)$$even if $\ann_{\D[s]}(f^s)\neq
\ann_{\D[s]}^{(1)}(f^s)$. For instance, when $X=\C^3$ and $f=x_1 x_2
(x_1+x_2) (x_1+x_2 x_3)$ (see example 6.2 in \cite{calde_nar_compo})
or in any of the examples in page 445 of \cite{cas_ucha_exper}. In
all this examples the divisor is free and satisfies the logarithmic
comparison theorem.
\end{ejemplo}

\section{Integrable logarithmic connections along quasi-homogeneous plane curves}
\label{sec:3}

Let $D\subset X=\C^2$ be a divisor defined by a reduced polynomial
equation $h(x_1,x_2)$, which is quasi-homogeneous with respect to
the strictly positive integer weights $\om_1,\om_2$ of the variables
$x_1,x_2$. We denote by $\om(f)$ the weight of a quasi-homogeneous
polynomial $f(x_1,x_2)$. The divisor $D$ is free, a global basis of
$\Der(\log D)$ is $\{\delta_1, \delta_2\}$, where
$$\textstyle\left( \begin{array}{c} \delta_1 \\
\delta_2
\end{array} \right) =
\left( \begin{array}{ccc}
 \om_1 x_1 & \om_2 x_2 \\
 - h_{x_2} &  h_{x_1}
\end{array} \right)
\left( \begin{array}{c} \dx{1} \\ \dx{2}
\end{array} \right).$$
We have:\smallskip

 \noindent -) $\delta_1(h)= \omega(h) h,\quad
\delta_2(h)=0,$\\
-) the determinant of the coefficient matrix is equal to $\omega(h) h$,\\
-) $ [\delta_1,\delta_2] = c \delta_2$, with $c = \omega(h) - \om_1 -
\om_2.$
\medskip

 We consider a logarithmic
connection ${\E}=\oplus_{i=1}^n\OX e_i$ given by actions:
$$\delta_1\cdot \left( \begin{array}{c} e_1 \\ \vdots \\
e_n
\end{array} \right) =
A_1 \left( \begin{array}{c} e_1 \\ \vdots \\
e_n
\end{array} \right), \quad
\delta_2 \cdot \left( \begin{array}{c} e_1 \\ \vdots \\
e_n
\end{array} \right) =
A_2 \left( \begin{array}{c} e_1 \\ \vdots  \\ e_n
\end{array} \right).
$$
For ${\E}$ to be integrable, the following integrability condition
\begin{equation} \label{eq:integrab1}\delta_1(A_2)-\delta_2(A_1) + [A_2,A_1]=c A_2
\end{equation}
must hold.
\medskip

\numero \label{nume:focus} We shall focus on the case where $A_1,A_2$
are $n\times n$ matrices satisfying (\ref{eq:integrab1}) and of the
form: {\small$$ A_1  = \left( \begin{array}{cccccc}
 -a &  0 & 0 & \cdots & 0 & 0\\
 -\delta_2(a) & -a\!+\!c  & 0&\cdots  &0&0\\
 -\delta_2^2(a) & -\comb{2}{1} \delta_2(a) & -a\!+\!2c &
 \cdots & 0 & 0\\
 \vdots & \vdots & \vdots & \ddots &\vdots&\vdots \\
 -\delta_2^{n-2}(a) & -\comb{n\!-\!2}{1} \delta_2^{n-3}(a)& -\comb{n\!-\!2}{2}
 \delta_2^{n-4}(a) & \cdots & - a\!+\!(n\!-\!2)c& 0\\
-\delta_2^{n-1}(a) & -\comb{n\!-\!1}{1} \delta_2^{n-2}(a)&
-\comb{n\!-\!1}{2}
 \delta_2^{n-3}(a) & \cdots & -\comb{n\!-\!1}{n\!-\!2}
\delta_2(a)&  -a\!+\!(n\!-\!1)c
\end{array} \right),
$$}
$$
A_2  = \left( \begin{array}{ccccc}
 0 &  1 & 0 & \cdots & 0\\
 0 & 0 & 1 & \cdots & 0 \\
 \vdots & \vdots & \vdots & \ddots & 0 \\
0 & 0 & 0 & \cdots & 1 \\
 -b_0 & -b_1 & -b_2 &\cdots & -b_{n-1}
\end{array} \right).
$$
with $a,b_0,\dots,b_{n-1}$ polynomials. Let us call $\E_{a,\ub}$ the
corresponding ILC.
\smallskip

\begin{lema} \label{lema:ciclico}
The $\VCERO$-module ${\E_{a,\ub}}$ is generated by $e_1$ (so it is
cyclic) and the $\VCERO$-annihilator of $e_1$ is the left ideal
$J_{a,\ub}$ generated by $\delta_1+a$ and
$\delta_2^n+b_{n-1}\delta_2^{n-1}+\cdots + b_1\delta_2 + b_0$. So,
the $\VCERO$-module ${\E_{a,\ub}}$ is isomorphic to
$\VCERO/J_{a,\ub}$.
\end{lema}

\begin{prueba} The first part is clear since $\delta_2\cdot e_i =
e_{i+1}$ for $i=1,\dots, n-1$. For the second part, the inclusion $
J_{a,\ub}\subset\ann_{\VCERO}(e_1)$ is also clear. To prove the
opposite inclusion, we use the fact that any germ of logarithmic
differential operator $P$ has a unique expression as a sum
$P=\sum_{i,j}a_{i,j}\delta_1^{i}\delta_2^{j}$, where the $a_{i,j}$
are germs of holomorphic functions (\cite{calde_ens}, th.~2.1.4) and
a division argument.
\end{prueba}

\begin{nota} Theorem 2.1.4 in \cite{calde_ens} says that $\VCERO
= \OX[\delta_1,\delta_2]$ with relations:
$$ [\delta_1,f] = \delta_1(f), [\delta_2,f]=\delta_2(f), [\delta_1,\delta_2] = c
\delta_2,\quad f\in \OX.$$ In particular, we can define the {\em
support} and the {\em exponent} of any germ of logarithmic
differential operator $P$ (or of any polynomial logarithmic
differential operator in the Weyl algebra) by using the (unique)
expression $P=\sum_{i,j}a_{i,j}\delta_1^{i}\delta_2^{j}$, and we
obtain a division theorem and a notion of {\em Gr\"obner basis} for
ideals. Under this scope, the integrability condition
(\ref{eq:integrab1}) reads out as the fact that the generators
$$ g_1=\delta_1+a,\quad g_2=\delta_2^n+b_{n-1}\delta_2^{n-1}+\cdots +
b_0$$ of $J_{a,\ub}$ satisfy Buchberger's criterion, i.e. that
$\delta_2^n g_1 - \delta_1 g_2$ has a vanishing remainder with
respect to the division by $g_1, g_2$, and then they form a Gr\"obner
basis of $J_{a,\ub}$.
\end{nota}

\begin{corolario} \label{DIa}
The $\DX$-module $\DX\otimes_{\VCERO}{\E_{a,\ub}}$ is isomorphic to
$\DX/I_{a,\ub}$, where $I_{a,\ub}=
\DX(\delta_1+a,\delta_2^n+b_{n-1}\delta_2^{n-1}+\cdots + b_0)$.
\end{corolario}
\medskip

\noindent For any integer $k$, we can consider the logarithmic
connections $\E_{a,\ub}(kD)$ and $\E_{a,\ub}^*$ (see section
\ref{nume:1}).

\begin{lema} \label{kD} With the above notations, the ILC $\E_{a,\ub}(kD)$ and
$\E_{a+\om(h) k,\ub}$ are isomorphic.
\end{lema}
\begin{prueba} An $\OX$-basis of ${\E_{a,\ub}}(kD)$ is $\{e_i^k= e_i \otimes h^{-k}\}_{i=1}^n$
and the action  of $\Der(\log D)$ over this basis is given by (see
(\ref{eq:oper-ext})):
$$\delta_1\cdot e_i^k = (\delta_1\cdot e_i) \otimes h^{-k} + e_i
\otimes (-\om(h) k h^{-k}), \quad \delta_2 \cdot e_i^k =
(\delta_2\cdot e_i) \otimes h^{-k}.$$ Then, the isomorphism of
$\OX$-modules
$$\sum_{i=1}^n b_i e_i \in \E_{a+\om(h) k,\ub} \mapsto  \sum_{i=1}^n b_i e_i^k \in
 {\E_{a,\ub}}(kD)$$ is clearly
$\VCERO$-linear.
\end{prueba}
\medskip

\noindent The proof of the following proposition is clear.

\begin{proposicion} \label{PtoPf}
The morphism
$$ \varrho_{\E_{a,\ub},k,k'}: \DX\otimes_{\VCERO} \E_{a,\ub}((1-k')D) \to \DX\otimes_{\VCERO}
\E_{a,\ub}(kD), $$ defined in (\ref{eq:2}), corresponds, through the
isomorphisms in corollary \ref{DIa} and lemma \ref{kD}, to the
morphism
$$\varrho_{\E_{a,\ub},k,k'}': \overline{P}\in \DX/I_{a+\om(h)(1-k'),\ub} \mapsto
\overline{Ph^{k+k'-1}}\in \DX/I_{a+\om(h)k}.$$
\end{proposicion}
\medskip

For the dual connection $\E_{a,\ub}^*$, in order to simplify, let us
concentrate on case $n=2$, where the integrability condition
(\ref{eq:integrab1}) reduces to:
\begin{equation} \label{eq:integrab2}
(\delta_1 - c)(b_1)=2\delta_2(a),\quad (\delta_1 - 2c)(b_0)=
\delta_2^2(a) + b_1\delta_2(a).
\end{equation}

\begin{lema} \label{dual} With the above notations, the ILC $\E_{a,\ub}^*$ and
${\E_{c-a,\ub^*}}$, with $\ub=(b_1,b_0)$ and
$\ub^*=(-b_1,b_0-\delta_2(b_1))$, are isomorphic.
\end{lema}
\begin{prueba} The action of $\Der(\log D)$
 over the dual basis
$\{e_1^*,e_2^*\}$ in $\E_{a,\ub}^*$ is given by:
$$(\delta_i\cdot e_j^*)(e_k) = \delta_i (e_j^*(e_k)) - e_j^* (\delta_i
e_k)= -e_j^* (\delta_i e_k),$$ for $i=1,2$ and $j,k=1,2$ (see
(\ref{eq:oper-ext})). Then
$$\delta_1 \left( \begin{array}{c} e_1^* \\  e_2^*
\end{array} \right) = -A_1^t
\left( \begin{array}{c} e_1^* \\ e_2^*
\end{array} \right), \quad \delta_2 \left( \begin{array}{c} e_1^*
\\ e_2^*
\end{array} \right) = -A_2^t
\left( \begin{array}{c} e_1^* \\ e_2^*
\end{array} \right).$$
 Choosing the new basis $\{w_1=e_2^*,w_2=-e_1^*+b_1 e_2^*\}$ of $\E_{a,\ub}^*$, we obtain
$$\delta_1 \left( \begin{array}{c} w_1 \\
w_2
\end{array} \right) =\cdots=  \left( \begin{array}{cc}
 a-c & 0 \\
 \delta_2(a) & a
\end{array} \right)
\left( \begin{array}{c} w_1 \\  w_2
\end{array} \right),$$ $$ \delta_2 \left( \begin{array}{c} w_1 \\
w_2
\end{array} \right) =\cdots = \left( \begin{array}{cc}
0&1\\ \delta_2(b_1)-b_0 & b_1\end{array} \right) \left(
\begin{array}{c} w_1 \\ w_2
\end{array} \right)$$
 and the isomorphism of $\OX$-modules
 $$\sum_{i=1}^2 b_i w_i \in \E_{a,\ub}^*\mapsto \sum_{i=1}^2 b_i e_i \in
\E_{c- a,\ub^*}$$ is clearly $\VCERO$-linear.
\end{prueba}

\section{Some explicit examples} \label{sec:4}

In this section we consider the case where $D\subset X=\C^2$ is
defined by the reduced equation $h=x_1^2-x_2^3$, and then
$\om(x_1)=3$, $\om(x_2)=2$, $\om(h)=6$ and the basis of $\Der(\log
D)$ is $\{\delta_1, \delta_2\}$, with
$$\textstyle\left( \begin{array}{c} \delta_1 \\
\delta_2
\end{array} \right) =
\left( \begin{array}{ccc}
 3x_1 & 2x_2 \\
 3x_2^2 & 2x_1
\end{array} \right)
\left( \begin{array}{c} \dx{1} \\ \dx{2}
\end{array} \right),$$
\noindent -) $\delta_1(h)= 6h,\quad
\delta_2(h)=0,$\\
-) the determinant of the coefficient matrix
is equal to $6h$,\\
-) $ [\delta_1,\delta_2]=\delta_2$ ($c=1$).
\medskip

\numero \label{nume:global} Since the ILC $\E_{a,\ub}$ and the ideals
$I_{a,\ub}$ in corollary \ref{DIa} are defined globally by
differential operators with polynomial coefficients and $D$ has a
global polynomial equation, the study of morphism
$$ \rho_{\E_{a,\ub},k}:\DX\Lotimes_{\VCERO} \E_{a,\ub}(kD) \to \E_{a,\ub}(\star D)$$
 can be done
globally at the level of the Weyl algebra $\W_2 =
\CC[x_1,x_2,\dx{1},\dx{2}]$.
\medskip

The integrability conditions in (\ref{eq:integrab2}) (for $n=2$)
become in our case
\begin{equation} \label{eq:integrab3}
(\delta_1 - 1)(b_1)=2\delta_2(a),\quad (\delta_1 - 2)(b_0)=
\delta_2^2(a) + b_1\delta_2(a).
\end{equation}
Once $a$ is fixed, it allows us to determine, uniquely, $b_1$  (the
operator $\delta_1 - 1$ is injective), and to also determine $b_0$ up
to a term $e x_2$, $e\in \CC$ (the kernel of the operator $\delta_1 -
2$ is generated by $x_2$). In order to simplify, let us take
$$a=\la+\m x_1+\n x_2,$$where $\umu= (\la,\m,\n)$ are complex
parameters, and then
$$ b_1=2\m x_2^2+2\n x_1$$
and
$$b_0=e x_2+3\n x_2^2+4\m x_1 x_2+\n^2 x_1^2+2\m \n x_1
x_2^2+\m^2 x_2^4,$$ with $e$ another complex parameter. For
convenience (see the rational factorization of $B(s)$ below), let us
consider another complex parameter $\f$ and make $e=\f-\f^2$.
\smallskip

Let us define the family of ILC of rank two, $\bE_{\f,\umu} :=
\E_{a,\ub}$ (see \ref{nume:focus}), with $a,b_0,b_1$ as above. We
have $\bE_{\f,\umu} = \VCERO\cdot e_1$ and $\ann_{\VCERO} e_1 =
\VCERO(g_1,g_2)$, with $g_1= \delta_1+a$ and $g_2=\delta_2^2 + b_1
\delta_2 + b_0$ (see lemma \ref{lema:ciclico}). It is clear that
$\bE_{\f,\umu}= \bE_{1-\f,\umu}$.
\medskip

The conclusion of corollary \ref{cor:annDefs} can be globalized and
we obtain $$\ann_{\DX[s]} (e_1 h^s) = \DX[s](\Phi(g_1),\Phi(g_2)) =
\DX[s](\delta_1+a-6s,g_2)$$and
$$\ann_{\W_2[s]} (e_1 h^s) = \W_2[s] (\delta_1+a-6s,g_2).$$

Let us consider the Weyl algebra with parameters
$$\textstyle \W' = \C\left[\la,\m,\n,\f,x_1,x_2,\dx{1},\dx{2}\right][s]$$
and the left ideal $I$ generated by
$$h,\quad \delta_1+a-6s,\quad \delta_2^2+b_1
\delta_2+b_0.$$ By a Gr\"obner basis computation with an elimination
order, for example, with the help of \cite{D-mod-M2}, we compute the
generator $B(s)$ of the ideal $I\cap \C[s]$ and operators
$P(s),C(s),D(s) \in \W'$ such that
$$B(s) = P(s) h + C(s)(\delta_1+a-6s) + D(s)(\delta_2^2+b_1
\delta_2+b_0).$$ We find
$$ B(s)=
\left(s-\frac{\la-5}{6}\right) \left(s-\frac{\la-8}{6}\right)
\left(s-\frac{\la-\f-6}{6}\right)
\left(s-\frac{\la+\f-7}{6}\right).$$

For $\la,\f\in \C$, let us call $B_{\la,\f}(s)\in \C[s]$ the
polynomial obtained from $B(s)$ in the obvious way. We obtain then
for each $\f,\la,\m,\n\in \C$ the global Bernstein-Sato functional
equation
\begin{equation}\label{eq:BS} B_{\la,\f}(s) e_1 h^s  = P(s) \left( e_1
h^{s+1} \right)\end{equation} in $\bE_{\f,\umu}[h^{-1},s]h^s$.
Therefore, $b_{\bE_{\f,\umu},p}(s)\ |\ B_{\la,\f}(s)$ (see prop.
\ref{prop:BS}) for any $p\in D$\footnote{In fact it is possible to
show that $b_{\bE_{\f,\umu},0}(s)\ = B_{\la,\f}(s)$.} and
$$ \kappa(\bE_{\f,\umu}) \geq \tau(\la,\f):= \min \{\text{integer
roots of $B_{\la,\f}(s)$}\} \in \ZZ \cup \{+\infty\}.$$ We can apply
theorem \ref{teo:main} to deduce that morphism
$$ \rho_{\bE_{\f,\umu},k}: \DX\otimes_{\VCERO} \bE_{\f,\umu}(kD)
\xrightarrow{} \bE_{\f,\umu}(\star D)$$ is an isomorphism for all
$k\geq -\tau(\la,\f)$. On the other hand, from lemma \ref{dual} we
know that $ (\bE_{\f,\la,\m,\n})^* = \bE_{\f,1-\la,-\m,-\n}$ and then
morphism
$$ \rho_{\bE^*_{\f,\umu},k'}: \DX\otimes_{\VCERO}
\bE^*_{\f,\umu}(k'D) \xrightarrow{} \bE^*_{\f,\umu}(\star D)$$ is an
isomorphism for all $k'\geq -\tau(1-\la,\f)$.
\medskip

The above results can be rephrased in the following way:
\medskip

\noindent 1) Morphism
$$ \rho_{\bE_{\f,\umu},k}: \DX\otimes_{\VCERO} \bE_{\f,\umu}(kD)
\xrightarrow{} \bE_{\f,\umu}(\star D)$$ is an isomorphism if the four
following conditions hold:\smallskip

\noindent
$\la + 6k \neq -1,-7,-13,-19,\dots $\\
$\la +6k \neq 2,-4,-10,-16,\dots$ \\
$\la+6k-\f \neq 0,-6,-12,-18,\dots$\\
$\la+6k+\f \neq 1,-5,-11,-17,\dots$ \medskip

\noindent 2) Morphism
$$ \rho_{\bE^*_{\f,\umu},k'}: \DX\otimes_{\VCERO} \bE^*_{\f,\umu}(k'D)
\xrightarrow{} \bE^*_{\f,\umu}(\star D)$$ is an isomorphism if the
four following conditions hold:\smallskip

\noindent
$1-\la + 6k' \neq -1,-7,-13,-19,\dots $ \\
$1-\la +6k' \neq 2,-4,-10,-16,\dots $\\
$1-\la+6k'-\f \neq 0,-6,-12,-18,\dots $\\
$1-\la+6k'+\f \neq 1,-5,-11,-17,\dots $\smallskip

\noindent or equivalently, if the four following conditions
hold:\smallskip

\noindent
$\la-6k'\neq 2,8,14,20,\dots $\\
$\la - 6k' \neq -1, 5,11,17,\dots $\\
$\la+\f -6k' \neq 1,7,13,19,\dots $\\
$\la-\f-6k'\neq 1,-5,-11,-17,\dots $
\medskip

\noindent In particular, if the four following conditions:
\begin{enumerate}
\item[(i)]
$ \la \not\equiv 2\  (\Mod 6)$ or
 $\la=2$
\item[(ii)]
$ \la \not\equiv 5\ (\Mod 6)$ or $\la=-1$
\item[(iii)]
$ \la+\f \not\equiv 1\ (\Mod 6)$ or $\la+\f=1$
\item[(iv)]
$ \la-\f \not\equiv 0\ (\Mod 6)$ or $\la-\f=0$
\end{enumerate}
hold, both morphisms
$$ \rho_{\bE_{\f,\umu},1}: \DX\otimes_{\VCERO} \bE_{\f,\umu}(D)
\xrightarrow{} \bE_{\f,\umu}(\star D),$$
$$ \rho_{\bE^*_{\f,\umu},1}: \DX\otimes_{\VCERO} \bE^*_{\f,\umu}(D)
\xrightarrow{} \bE^*_{\f,\umu}(\star D)$$are isomorphisms.
\medskip

Let us denote by $\LL_{\f,\umu}$ the local system over
 $X-D$ of the horizontal sections of  $\bE_{\f,\umu}$. By theorem
 \ref{teo:2}, we have
  $$\IC_X ( \LL_{\f,\umu}) \simeq \DR(\im \varrho_{\bE_{\f,\umu},1,1}),$$
provided that conditions (i)-(iv) are satisfied.

Proposition \ref{PtoPf} and \ref{nume:global}
 reduce the computation of
$\im \varrho_{\bE_{\f,\umu},1,1}$ to the computation of the image of
the map
$$ \theta_{\f,\umu}:\overline{L}\in \W_2/\W_2(g_1,g_2) \mapsto
\overline{Lh}\in \W_2/\W_2(g_1+6,g_2),$$ but $\im \theta_{\f,\umu} =
\W_2/K_{\f,\umu}$ where
$$ K_{\f,\umu} = \{R\in \W_2\ |\ Rh \in \W_2(g_1+6,g_2)\}.$$

Now, in order to compute generators of $K_{\f,\umu}$, we proceed as
follows. Since $[g_1,g_2]=2g_2$ (for any $\f,\umu$) and the symbols
$\sigma(g_1)=\sigma(\delta_1)$, $\sigma(g_2)=\sigma(\delta_2)^2$ form
a regular sequence ($D$ is Koszul free!), we deduce that
$$ \sigma\left( \W_2(g_1+6,g_2^{\vphantom{1em}}) \right) =
\left(\sigma(\delta_1),\sigma(\delta_2)^2\right)$$ and consequently $
\sigma\left(K_{\f,\umu}\right) \subset
(\sigma(\delta_1),\sigma(\delta_2)^2): h$. A straightforward
(commutative) computation shows that
$$ (\sigma(\delta_1),\sigma(\delta_2)^2): h =
(\sigma(\delta_1),\sigma(Q_0))$$ with
$Q_0=9x_2\frac{\partial^2}{\partial
x_1^2}-4\frac{\partial^2}{\partial x_2^2}$, and
\begin{equation} \label{eq:clave}
\sigma(Q_0)h= x_2\sigma(\delta_1)^2 -
\sigma(\delta_2)^2=x_2\sigma(\delta_1)\sigma(g_1+6) - \sigma(g_2).
\end{equation}
Searching to lift the relation (\ref{eq:clave}) to $\W_2$, we find
$$ Qh = x_2(\delta_1+\m x_1+\n x_2+7-\la) (g_1+6)-g_2 +
(\la^2-\la+\f-\f^2)x_2,$$ with $ Q = Q_0 + 6\m x_2
\dx{1}-4\n\dx{2}+\m^2 x_2-\n^2$. In particular, if condition
\begin{equation}\label{eq:condition}
\la^2-\la+\f-\f^2 = 0\quad (\Leftrightarrow  \la-\f=0\ \ \text{or}\ \
\la+\f=1)
\end{equation}
holds, then $ Q\in K_{\f,\umu}$.

Actually, by using the equality $[Q,g_1]= 4Q$ and the fact that
$\sigma(Q)=\sigma(Q_0)$ and $\sigma(g_1)=\sigma(\delta_1)$ also form
a regular sequence in $\Gr \W_2$, condition (\ref{eq:condition})
implies that $$ K_{\f,\umu} = \W_2(g_1,Q),\quad
\sigma\left(K_{\f,\umu}\right) = (\sigma(\delta_1),\sigma(Q_0)).$$ On
the other hand, since $\sigma(Q_0)$ is not contained in the ideal
$(x_1,x_2)$, we finally deduce the following result:
\medskip

If parameters $\f,\umu=(\la,\m,\n)$ satisfy conditons (i)-(iv) and
(\ref{eq:condition}), then the conormal of the origin $T_0^*(X)$ does
not appear as an irreducible component of the characteristic variety
of $\im \theta_{\f,\umu} = \W_2/K_{\f,\umu}$, and consequently
$$\Ch (\IC_X(\LL_{\f,\umu})) = \Ch \left(\W_2/K_{\f,\umu}\right)= \{ \sigma({\delta_1})=
\sigma(Q_0)=0\}= T_X^*(X) \cup T_D^*(X).$$ The existence of such an
example has been suggested by \cite{nar_88}, example (3.4), but the
question on the values of the parameters $\f,\umu$ for which the
local system $\LL_{\f,\umu}$ is irreducible will be treated
elsewhere.
\medskip

If condition (\ref{eq:condition}) does not hold, it is not clear that
there exists a general expression for a system of generators of
$K_{\f,\umu}$ as before.

\begin{nota} The relationship between the preceding results and
examples and the hypergeometric local systems (cf.
\cite{neto_silva_cras_02,neto_silva_pacific_02,psilva_tesis}) is
interesting and possibly deserves further work.
\end{nota}

\def\cprime{$'$}

\bigskip

{\small \noindent Departamento de \'{A}lgebra\\
 Facultad de  Matem\'{a}ticas\\ Universidad de Sevilla\\ P.O. Box 1160\\ 41080
 Sevilla\\ Spain}. \\
{\small {\it E-mail}:  $\{$calderon,narvaez$\}$@algebra.us.es
 }
\end{document}